\documentclass[reqno]{amsart}

\usepackage{article-preamble}
\usetikzlibrary{arrows.meta}

\hypersetup{colorlinks=true,
    linkcolor=blue,
    citecolor=blue,
    pdftitle={Bounds on the dimension of lineal extensions},
    pdfauthor={Ryan E.~G.~Bushling and Jacob B.~Fiedler}}

\begin{document}

\title[Bounds on the dimension of lineal extensions]{Bounds on the dimension of lineal extensions}
\author{Ryan E.~G.~Bushling}
\address{Department of Mathematics \\ University of Washington, Box 354350 \\ Seattle, WA 98195-4350}
\email{reb28@uw.edu}
\author{Jacob B.~Fiedler}
\address{Department of Mathematics \\ University of Wisconsin, 480 Lincoln Drive \\ Madison, WI 53706-1325}
\email{jbfiedler2@wisc.edu}
\thanks{The second author was supported in part by NSF DMS-2037851 and NSF DMS-2246906.}

\begin{abstract}
    Let $E \subseteq \mathbb{R}^n$ be a union of line segments and $F \subseteq \mathbb{R}^n$ the set obtained from $E$ by extending each line segment in $E$ to a full line. Keleti's \textit{line segment extension conjecture} posits that the Hausdorff dimension of $F$ should equal that of $E$. Working in $\mathbb{R}^2$, we use effective methods to prove a strong packing dimension variant of this conjecture. Furthermore, a key inequality in this proof readily entails the planar case of the generalized Kakeya conjecture for packing dimension. This is followed by several doubling estimates in higher dimensions and connections to related problems.
\end{abstract}

\maketitle

\section{Introduction and main results}

Let $E = \bigcup \mathcal{I}$, where $\mathcal{I}$ is a family of line segments in $\R^n$, $n \geq 2$. Throughout, $\mathcal{I}$ is assumed to be \textit{maximal} in the sense that, if $I$ is a line segment and $I \subseteq \bigcup \mathcal{I}$, then $I \in \mathcal{I}$\textemdash a hypothesis that results in no loss of generality in what follows. Denoting by $\mathcal{A}(n,1)$ the affine Grassmannian of lines in $\R^n$, we define the \define{lineal extension of $\bm{E}$} to be the set $\upbold{L}(E)$ formed from $E$ by extending each $I \in \mathcal{I}$ to the unique line $\ell_I \subset \R^n$ containing $I$:
\begin{equation*}
    \upbold{L}(E) \ := \ \bigcup_{I \in \mathcal{I}} \ell_I \ = \ \bigcup \big\{ \ell \in \mathcal{A}(n,1) \!: E \cap \ell \text{ contains a line segment} \big\}.
\end{equation*}
With this setup, Keleti \cite{keleti2016lines} proposed the following conjecture. Let $\hdim$ denote Hausdorff dimension.

\begin{conj}[Line segment extension conjecture] \label{conj:lse}
    Let $E \subseteq \R^n$ be a union of line segments and $\upbold{L}(E)$ its lineal extension. Then $\hdim \upbold{L}(E) = \hdim E$.
\end{conj}

Conjecture \ref{conj:lse} is open in dimensions $n \geq 3$ but is known for $n=2$.

\begin{thm}[Keleti \cite{keleti2016lines}] \label{thm:lse-R2}
    If $E \subseteq \R^2$ is a union of line segments, then $\hdim \upbold{L}(E) = \hdim E$.
\end{thm}

This article concerns variants and extensions of this problem, emphasizing the planar case.

\subsection{Line segment extension in \texorpdfstring{$\bm{\R^2}$}{R2}}

We begin by introducing some notation. Let $\pdim$ denote packing dimension.

\begin{defn} \label{defn:extensions}
    For $s \in (0,1]$ and $E \subseteq \R^n$ a union of subsets of lines with (packing or Hausdorff) dimension at least $s$, let
    \begin{align*}
        L_s^{\rm H}(E) &:= \bigcup \big\{ \ell \in \mathcal{A}(n,1) \!: \hdim (\ell \cap E) \geq s \big\}
        \intertext{and}
        L_s^{\rm P}(E) &:= \bigcup \big\{ \ell \in \mathcal{A}(n,1) \!: \pdim (\ell \cap E) \geq s \big\}
    \end{align*}
    be, respectively, the \define{$\bm{s}$-Hausdorff extension} and \define{$\bm{s}$-packing extension of $\bm{E}$}.\footnote{Alternatively, in the above definitions, we could drop the requirement that $E$ be the union of certain subsets of lines and instead write $L_s^{\rm H}(E) := \bigcup \big\{ \ell \in \mathcal{A}(n,1) \!: \hdim (E \cap \ell) \geq s \big\} \cup E$, and likewise for packing dimension. One could partition $E$ into (1) the set of points covered by lines intersecting $E$ in sets of Hausdorff dimension at least $s$ and (2) the remainder of $E$; calling the first of these $E_\text{lines}$, it is easy to see that $L_s^{\rm H}(E)$=$L_s^{\rm H} (E_\text{lines}) \cup E$. The same holds for packing dimension. Morally speaking, then, the actual definitions we use encompasses all of the interesting features of the problem.} In the extreme case $s = 0$ we let
    \begin{equation*}
        L_0(E) := \bigcup \big\{ \ell \in \mathcal{A}(n,1) \!: \# (\ell \cap E) \geq 2 \big\},
    \end{equation*}
    which we call the \define{two-point extension of $\bm{E}$}.
\end{defn}

Our first result\textemdash a generalization of Theorem \ref{thm:lse-R2}\textemdash pertains to the lineal extension of Furstenberg sets.

\begin{prop} \label{prop:generalized-lse-R2}
    Let $s \in (0,1]$ and let $E \subseteq \R^2$ be a union of (at least) $s$-Hausdorff-dimensional subsets of lines. Then $\hdim L_s^{\rm H}(E) \leq \hdim E + 2 - 2s$. In particular, if $E$ is a union of line segments and $\upbold{L}(E)$ is its lineal extension, then $\hdim \upbold{L}(E) = \hdim E$.
\end{prop}

Keleti's proof of Theorem \ref{thm:lse-R2} combines Marstrand's slicing theorem and the ``Fubini inequality" for Hausdorff measures in a simple and elegant argument, whereas we prove Proposition \ref{prop:generalized-lse-R2} by effective methods that, in particular, hinge on a result \cite{lutz2020bounding} of N.~Lutz and Stull (Theorem \ref{thm:bounding} below). We remark here that this is in fact implied by the Furstenberg set bound \cite{molter2012furstenberg} of Molter and Rela; cf.~\S \ref{s:lse-background} below.

One motivation for Proposition \ref{prop:generalized-lse-R2} is that the proof is morally similar to (but much simpler than) that of our main result.

\begin{thm} \label{thm:packing-lse-R2}
    If $E \subseteq \R^2$ is a union of $1$-Hausdorff-dimensional subsets of lines, then $\pdim L_1^{\rm H}(E) = \pdim E$. In particular, if $E$ is a union of line segments and $\upbold{L}(E)$ is its lineal extension, then $\pdim \upbold{L}(E) = \pdim E$.
\end{thm}

This theorem follows from a different effective analogue that cannot be proved directly from the aforementioned Lutz--Stull result. The bulk of the proof in \S \ref{s:plane} is establishing this analogue, which involves a sort of multiscale application of the ideas underlying \cite{lutz2020bounding}.

A strong ``generalized Kakeya conjecture" for packing dimension (see \S \ref{ss:history}) follows readily from the penultimate step of the proof of Theorem \ref{thm:packing-lse-R2}. Let $\bbp^{n-1} := \bbs^{n-1}/\{\pm 1\}$ be the set of directions of lines in $\R^n$.

\begin{cor} \label{cor:packing-gen-kakeya}
    Let $E \subseteq \R^2$ and let $D \subseteq \bbp^1$ be the set of directions of lines intersecting $E$ in a set of Hausdorff dimension $1$. If $D \neq \varnothing$, then
    \begin{equation} \label{eq:packing-gen-kakeya}
        \pdim D + 1 \leq \pdim E.
    \end{equation}
\end{cor}

\subsection{Line segment extension in \texorpdfstring{$\bm{\R^n}$}{Rn} and elementary Besicovitch set estimates}

Partial results in higher dimensions, including some (non-)doubling bounds on the dimension of lineal extensions, follow from more rudimentary ``two-part code" arguments in the spirit of \cite{altaf2023distance} Theorem 1.2.

\begin{prop} \label{prop:doubling-bound}
    If $E \subseteq \R^n$ is a union of line segments, then
    \begin{equation} \label{eq:doubling-bound}
        \hdim \upbold{L}(E) \leq \hdim E + \pdim E - 1 \quad \text{and} \quad \pdim \upbold{L}(E) \leq 2 \pdim E - 1.
    \end{equation}
\end{prop}

Such results are connected to the Kakeya conjecture via the following theorem. Call a subset of $\R^n$ (not necessarily Borel) a \define{Besicovitch set} if it contains a unit line segment in every direction.

\begin{thm}[Keleti \cite{keleti2016lines}] \label{thm:lse-to-kakeya} \hs{0}
\begin{enumerate}[label={\normalfont \textbf{\alph*.}}, noitemsep, topsep=-3pt]
    \item If the line segment extension conjecture holds in $\R^n$ for {\rm some} $n \geq 2$, then every Besicovitch set in $\R^n$ has Hausdorff dimension at least $n-1$.
    \item If the line segment extension conjecture holds in $\R^n$ for {\rm all} $n \geq 2$, then, for every $n \geq 2$, every Besicovitch set in $\R^n$ has packing dimension $n$.
\end{enumerate}
\end{thm}

While Theorem \ref{thm:lse-to-kakeya} assumes the full strength of Conjecture \ref{conj:lse}, with a small modification of the final step in Keleti's proof\textemdash to which the reader is referred\textemdash one obtains the following generalization.

\begin{lem} \label{lem:lse-to-kakeya}
    Suppose there is a function $g \!: [0,n]^2 \to [0,n]$ such that the following holds: \vs{-0.2}
    \begin{quote}
        {\rm If $E$ is a union of line segments in $\R^n$, then $g(\hdim E, \pdim E) \geq \hdim \upbold{L}(E)$.} \vs{-0.65}
    \end{quote}
    Then $g(\hdim K, \pdim K) \geq n-1$ for every Besicovitch set $K \subseteq \R^n$.
\end{lem}

As an immediate consequence of this lemma and Proposition \ref{prop:doubling-bound}, we obtain an elementary estimate on the dimension of Besicovitch sets in $\R^n$. (See also \S \ref{s:lse-background} for implications for the generalized Kakeya conjecture.)

\begin{cor} \label{cor:kakeya-bound}
    If $K \subseteq \R^n$ is a Besicovitch set, then
    \begin{equation*}
        \hdim K + \pdim K \geq n.
    \end{equation*}
\end{cor}

This is of course far from state-of-the-art, but the implication of Kakeya inequalities from line segment extension inequalities has something of a practical importance that we describe below.

\section{Background on the line segment extension conjecture and its relatives} \label{s:lse-background}

\subsection{History and context} \label{ss:history}

Keleti \cite{keleti2016lines} introduced the line segment extension problem as a natural follow-up to the constructions of Nikodym \cite{nikodym1927mesure} and Larman \cite{larman1971compact} showing that a union of closed line segments in $\R^n$ can have positive Lebesgue measure even when the union of the corresponding open line segments is Lebesgue-null.

In this same vein, Falconer and Mattila \cite{falconer2016strong} introduced a ``hyperplane extension problem," which they treated as a slicing problem with a dual projection problem that is amenable to a Marstrand-type exceptional set estimate. The planar case of their Theorem 3.2 weakens the hypothesis in Theorem \ref{thm:lse-R2} that $E$ contains many line segments to the hypothesis that it contains many positive-measure subsets of lines, which Proposition \ref{prop:generalized-lse-R2} further weakens to $s$-dimensional subsets of lines (possibly of $\mathcal{H}^s$-measure $0$). Another consequence of \cite{falconer2016strong} is an equation for the dimension of a family of hyperplanes in terms of that of its union. H\'{e}ra, Keleti, and M\'{a}th\'{e} \cite{hera2019hausdorff} pursued this direction in arbitrary dimension and codimension, bounding the dimension of families $\Lambda \subseteq \mathcal{A}(n,k)$ from above in terms of the dimension of any set giving large slices to $\bigcup \Lambda$.

Redirecting attention back to the connection between the line segment extension problem and the Kakeya problem, Keleti and M\'{a}th\'{e} \cite{keleti2022equivalences} showed that Theorem \ref{thm:lse-to-kakeya} has a strong converse.

\begin{thm}[Keleti--M\'{a}th\'{e} \cite{keleti2022equivalences}] \label{thm:kakeya-to-lse}
    If the Kakeya conjecture is true in $\R^n$, then the line segment extension conjecture is true in $\R^n$.
\end{thm}

This they established as a corollary to the equivalence of the Kakeya conjecture with the \textit{generalized Kakeya conjecture}.

\begin{conj}[Generalized Kakeya conjecture] \label{conj:gen-kakeya}
    Let $E \subseteq \R^n$ and let $D \subseteq \bbp^{n-1}$ be the set of directions in which $E$ contains a line segment. If $D \neq \varnothing$, then
    \begin{equation*}
        \hdim D + 1 \leq \hdim E.
    \end{equation*}
\end{conj}

In particular, this conjecture is true in $\R^2$ but open in all higher dimensions. Unlike Conjecture \ref{conj:gen-kakeya}, the generalized Kakeya conjecture for packing dimension does not seem to readily imply the packing dimension analogue of Conjecture \ref{conj:lse}. On the Hausdorff side this implication follows from Marstrand's slicing theorem, of which the packing dimension analogue is considerably weaker (cf.~\cite{falconer1996packing}). In fact, without additional hypotheses, replacing Hausdorff dimension with packing dimension in the proof of Conjecture \ref{conj:lse} from Conjecture \ref{conj:gen-kakeya} only gives the trivial lower bound $\pdim E \geq 1$. It is furthermore not obvious to us that the argument used to establish the equivalence between Kakeya and generalized Kakeya in \cite{keleti2022equivalences} easily adapts to packing dimension, and for these reasons it seems surprising that both Theorem \ref{thm:packing-lse-R2} and Corollary \ref{cor:packing-gen-kakeya} fall out of a single proof.

On a different note, it should also be remarked here that the implication in Theorem \ref{thm:kakeya-to-lse} is \textit{not} quantitative, in the sense that absolute lower bounds on the size of Besicovitch sets (or of unions of line segments more generally) do not translate into progress toward the line segment extension conjecture. This stands in contrast to Lemma \ref{lem:lse-to-kakeya}, according to which partial results on the line segment conjecture \textit{do} make headway on the Kakeya problem. In fact, the equivalence between the Kakeya and generalized Kakeya conjectures proved in \cite{keleti2022equivalences} is also quantitative:

\begin{thm}[Keleti--M\'{a}th\'{e} \cite{keleti2022equivalences}] \label{thm:kakeya-v-gen-kakeya}
    Let $E \subseteq \R^n$ and let $\varnothing \neq D \subseteq \bbp^{n-1}$ be the set of directions in which $E$ contains a line segment. Then there exists a compact Besicovitch set $K \subset \R^n$ with
    \begin{equation*}
        \hdim K \leq n-1 + \hdim E - \hdim D.  
    \end{equation*}
\end{thm}

Our method gives such inequalities of generalized Kakeya type for packing and mixed Hausdorff-packing dimensions directly, without reference to a general result analogous to Theorem \ref{thm:kakeya-v-gen-kakeya}; see the remarks following the proof of Proposition \ref{prop:doubling-bound} in \S \ref{s:higher-dim}.

\subsection{\texorpdfstring{$\bm{s}$}{s}-Hausdorff extensions and the Furstenberg set conjecture}

We take a moment to expound on the relationship between Proposition \ref{prop:generalized-lse-R2} and the Furstenberg set conjecture. For $s \in (0,1]$ and $t \in [0,2]$, call a set $E \subseteq \R^2$ an \define{$\bm{(s,t)}$-Furstenberg set} if there exists a nonempty $t$-Hausdorff-dimensional set $\Lambda \subseteq \mathcal{A}(2,1)$ such that
\begin{equation*}
    E = \bigcup_{\ell \in \Lambda} (E \cap \ell), \qquad \text{where} \qquad \hdim(E \cap \ell) \geq s \hs{0.25} \forall \+ \ell \in \Lambda.
\end{equation*}
If $\Lambda$ is the maximal set of lines with this property, then a special case of the aforementioned \cite{falconer2016strong} Theorem 3.2 gives
\begin{equation} \label{eq:falconer-mattila}
    \hdim L_s^{\rm H}(E) = \min \, \{t+1,2\},
\end{equation}
which in conjunction with Proposition \ref{prop:generalized-lse-R2} entails
\begin{equation*}
    \min \, \{ t+1, 2 \} \leq \hdim E + 2 - 2s, \qquad \text{i.e.,} \qquad 2s + \min \, \{t,1\} - 1 \leq \hdim E.
\end{equation*}
When $E$ is Borel, this is essentially the $1 \leq s + \tfrac{t}{2}$ case of the Molter and Rela \cite{molter2012furstenberg} Furstenberg set bound, and running the argument in reverse in turn yields Proposition \ref{prop:generalized-lse-R2} for Borel sets from their Furstenberg set estimate.\footnote{The authors thank Tam\'{a}s Keleti and Joshua Zahl for sharing this observation.}

More recently, Ren and Wang fully resolved the Furstenberg set conjecture in the plane.

\begin{thm}[Ren--Wang \cite{ren2023furstenberg}] \label{thm:ren-wang}
    If $E \subseteq \R^2$ is a Borel $(s,t)$-Furstenberg set, then
    \begin{equation*} \label{eq:ren-wang}
        \hdim E \geq \min \left\{ s+t, \frac{3s+t}{2}, s+1 \right\}.
    \end{equation*}
\end{thm}

A corollary, then, is a strengthening of Proposition \ref{prop:generalized-lse-R2} for Borel sets. The sharp examples for the Ren--Wang inequality are likewise sharp for this corollary.

\begin{cor} \label{cor:ren-wang}
    Let $s \in (0,1]$ and let $E \subseteq \R^2$ be a union of (at least) $s$-dimensional subsets of lines. If $E$ is Borel, then
    \begin{equation} \label{eq:ren-wang-cor}
        \hdim L_s^{\rm H}(E) \leq \max \big\{ \hdim E + 1 - s, \, 2 \hdim E + 1 - 3s \big\}.
    \end{equation}
\end{cor}

\textit{Proof.} We work by cases according to the values of $s$ and $t := \hdim \, \{ \ell \in \mathcal{A}(2,1) \!: \hdim(E \cap \ell) \geq s \}$. Suppose first that $s+t \geq 2$. Then necessarily $t \geq 1$, so Theorem \ref{thm:ren-wang} implies $\hdim E \geq s+1$ and \eqref{eq:falconer-mattila} implies $\hdim L_s^{\rm H}(E) = 2$, from which it follows that
\begin{equation*}
    \hdim L_s^{\rm H}(E) = (s+1) + 1 - s \leq \hdim E + 1 - s.
\end{equation*}
Next, suppose instead that $s+t \leq 2$ and $s \leq t$. Then Theorem \ref{thm:ren-wang} implies $\hdim E \geq \tfrac{3s+t}{2}$ and \eqref{eq:falconer-mattila} gives $\hdim L_s^{\rm H}(E) \leq t+1$, so
\begin{equation*}
    \frac{3s + \hdim L_s^{\rm H}(E) - 1}{2} \leq \frac{3s+t}{2} \leq \hdim E.
\end{equation*}
Rearranging terms gives the second expression on the right-hand side of \eqref{eq:ren-wang-cor}.

Finally, suppose that $s+t \leq 2$ and $0 \leq t < s$, so that $\hdim E \geq s+t$ by Theorem \ref{thm:ren-wang} and $\hdim L_s^{\rm H}(E) = t+1$ by \eqref{eq:falconer-mattila}. Then
\begin{equation} \label{eq:t-leq-s}
    s + \hdim L_s^{\rm H}(E) - 1 = s+t \leq \hdim E,
\end{equation}
once more giving the first expression on the right-hand side of \eqref{eq:ren-wang-cor}. As this covers all three cases, the corollary is proved. \textqed

Conversely, the final step in the proof of Proposition \ref{prop:generalized-lse-R2} shows that the conclusion \eqref{eq:t-leq-s} holds when $s \geq t$, recovering the corresponding case of Theorem \ref{thm:ren-wang}.

\section{Preliminaries on effective methods} \label{s:complexity}

\subsection{Basic definitions} The main goal of this section is to collect in one place several tools which we use repeatedly in the remainder of the paper, especially for the benefit of readers less familiar with Kolmogorov complexity. We operate in the algorithmic framework laid out in \cite{lutz2018algorithmic}, which we briefly recall here to establish terminology and notation. Let $\{0,1\}^*$ be the collection of all finite strings over $\{0,1\}$, including the empty string $\emptyset$. Fixing some prefix-free universal oracle Turing machine $U$, we define for each pair $\sigma,\tau \in \{0,1\}^*$ the \define{Kolmogorov complexity of $\bm{\sigma}$ given $\bm{\tau}$} to be the minimal length of a program that, when given to $U$ as an input with side information $\tau$, returns $\sigma$ as the output:
\begin{equation*}
    K(\sigma \vert \tau) := \min \big\{ |\pi| \in \N \!: \pi \in \{0,1\}^*, \, U(\pi,\tau) = \sigma \big\}.
\end{equation*}
When $\tau = \emptyset$, we write $K(\sigma) := K(\sigma\vert\emptyset)$ and simply call this quantity the \define{Kolmogorov complexity of $\bm{\sigma}$}.

The ``universality" of $U$ refers to the fact that, for every prefix-free oracle Turing machine $M$, there exists a program $\pi_M \in \{0,1\}^*$ such that
\begin{equation*}
    U(\pi_M,\sigma) = M(\sigma) \qquad \forall \+ \sigma \in \{0,1\}^*.
\end{equation*}
The length of the shortest such $\pi_M$ is called the \textit{machine constant of $M$}. When it is more awkward to work with $U(\pi_M, \, \cdot \,)$ than it is to work with $M$ directly, we opt for the latter and then add the machine constant to the length of the shortest $\sigma$ such that $M(\sigma) = \tau$ when computing the Kolmogorov complexity of $\tau$.

Identifying the family of all rational vectors with $\{0,1\}^*$ via some \textit{effective encoding} $\bigcup_{n \in \N} \Q^n$ $\hookrightarrow \{0,1\}^*$, we may extend these definitions from strings to real vectors as follows. Let $x \in \R^n$, $y \in \R^m$, and $r,s \in \N$. We call
\begin{equation*}
    K_r(x) \, := \min_{p \in B(x,2^{-r}) \cap \Q^n} K(p)
\end{equation*}
the \define{Kolmogorov complexity of $\bm{x}$ to precision $\bm{r}$} and
\begin{equation} \label{eq:conditional-complexity}
    K_{r,s}(x \vert y) \, := \max_{q \in B(y,2^{-s}) \cap \Q^m} \left( \min_{p \in B(x,2^{-r}) \cap \Q^n} K(p \vert q) \right)
\end{equation}
the \define{Kolmogorov complexity of $\bm{x}$ to precision $\bm{r}$ given $\bm{y}$ to precision $\bm{s}$}. When $s=r$ we simply write $K_r(x \vert y) := K_{r,r}(x \vert y)$, and when $y=x$ we write $K_{r,s}(x) := K_{r,s}(x \vert x)$.

When working with conditional complexities, we will have occasion to make statements like, ``$\pi$ testifies to $K_{r,s}(x \vert y)$." By this we mean, ``$\pi$ testifies to $K(x \-\upharpoonright\- r \vert y \-\upharpoonright\- s)$," where $\upharpoonright$ denotes the truncation of a (possibly infinite) string to the given number of bits of precision. This technical nuance is necessary owing to the fact that there may be many pairs $(p,q)$ as in \eqref{eq:conditional-complexity} such that $K_{r,s}(x \vert y) = K(p \vert q)$. In particular, there is generally no single $p \in (x,2^{-r}) \cap \Q^n$ such that $K(p \vert q) \leq K_{r,s}(x \vert y)$ for \textit{all} $q \in B(y,2^{-s}) \cap \Q^m$. However, since $K_{r,s}(x \vert y) = K(x \-\upharpoonright\- r \vert y \-\upharpoonright\- s) + O(\log r + \log s)$ (cf.\! \cite{lutz2020bounding} Corollary 2.5), one can work as though such $p$ did exist, modulo a logarithmic error.

By allowing a machine access to an \textit{oracle} $A \subseteq \{0,1\}^*$, we can \textit{relativize} the above definitions to $A$, in which case we embellish the symbols $U$, $K$, and $M$ with a superscript $A$. An oracle represents extra information that an oracle Turing machine is allowed to use in computations. Access to an oracle can never make a computation meaningfully harder, as a machine can always ``ignore" the oracle if its information is irrelevant. In particular, if $A$ and $B$ are oracles, then
\begin{equation*}
    K_r^{A,B}(x) \leq K_r^A(x) + O(1)
\end{equation*}
for all $x \in \R^n$.

Using some standard encoding, we can consider points in $\R^n$ as oracles. Intuitively, conditional access to a point up to a certain precision should be no more useful than oracle access to \textit{all} of the information in that point, and this is made precise by the inequality
\begin{equation*}
    K_r^{A,x}(y) \leq K_r^A(y \vert x) + O(\log r).
\end{equation*}

\subsection{Some useful results} One key property of Kolmogorov complexity is \textit{symmetry of information}. The following quantitative form will see repeated use in this paper.

\begin{lem}[Symmetry of information \cite{lutz2020bounding}] \label{lem:symmetry}
    For all $A \subseteq \{0,1\}^*$, $x \in \R^n$, $y \in \R^m$, and $r,s \in \N$ with $r \geq s$: \vs{-0.15}
    \begin{enumerate}[label={\normalfont\textbf{\alph*.}},itemsep=3pt,topsep=0pt]
        \item $\big| K_r^A(x \vert y) + K_r^A(y) - K_r^A(x,y) \big| \leq O(\log r) + O(\log \log \|y\|)$.
        \item $\big| K_{r,s}^A(x) + K_s^A(x) - K_r^A(x) \big| \leq O(\log r) + O(\log \log \|x\|)$.
    \end{enumerate}
\end{lem}

In practice, the norms of the points we work with are fixed and independent of the level of precision, so we frequently use these facts in the (relativized) forms
\begin{equation*}
   K_r^A(x,y) \approx K_r^A(x \vert y) + K_r^A(y) \quad \text{and} \quad K_r^A(y) \approx K_{r,s}^A(y) + K_s^A(y),
\end{equation*}
where both equalities hold up to a logarithmic term in $r$. The latter of these is particularly useful as a tool to bound the complexity of $y$ at a given precision: its repeated use allows us to \textit{partition} the interval $[1,r]$ into smaller intervals on which the complexity function of $y$ may have more desirable properties.

We add to this another result for understanding how complexity varies with precision. Case and J.~Lutz \cite{case2015dimension} showed that, for any $A \subseteq \{0,1\}^*$, $r,s \in \N$, and $x \in \R^n$,
\begin{equation*}
    K_r^A(x) \leq K_{r+s}^A(x) \leq K_r^A(x) + ns + O(\log(s + r)).
\end{equation*}
This bound captures two essential features of the Kolmogorov complexity of points: it is non-decreasing, and its growth rate is essentially bounded by $n$ on sufficiently long intervals. 

Ultimately, we study the Kolmogorov complexity of points in $x \in \R^n$ to bound their asymptotic information density. Given an oracle $A \subseteq \{0,1\}^*$, we call
\begin{equation*}
    \dim^A(x) := \liminf_{r \to \infty} \frac{K_r^A(x)}{r} \quad \text{and} \quad \Dim^A(x) := \limsup_{r \to \infty} \frac{K_r^A(x)}{r}
\end{equation*}
the \define{effective Hausdorff dimension} and the \define{effective packing dimension of $\bm{x}$ relative to $\bm{A}$}, respectively. The utility of effective dimensions in geometric measure theory stems from the following theorem of J.~Lutz and N.~Lutz.

\begin{thm}[Point-to-set principle \cite{lutz2018algorithmic}] \label{thm:pts}
    For every set $E \subseteq \R^n$,
    \begin{equation*}
        \hdim E \, = \min_{A \subseteq \{0,1\}^*} \sup_{x \in E} \dim^A(x) \quad \text{ and } \quad \pdim E \, = \min_{A \subseteq \{0,1\}^*} \sup_{x \in E} \Dim^A(x).
    \end{equation*}
\end{thm}

We frequently use the following immediate consequence of this theorem: given some $E \subseteq \R^n$, for any oracle $A$ and $\eps > 0$, there exists some $x \in E$ such that $\dim^A(x) > \hdim E - \eps$, and likewise for packing dimension. 

Turning our attention to the effective dimension of points on lines, we note the following observation of N.~Lutz and Stull: for any $A \subseteq \{0,1\}^*$ and any $x,a,b \in \R$, 
\begin{equation*}
    K_r^A(x, ax+b) \leq K_r^A(x,a,b) + O_{x,a,b}(\log r).
\end{equation*}
This is because, for every large enough precision $r$, a Turing machine given $x,a,b$ at precision $r$ can perform a very accurate calculation of $ax+b$ at precision $r$. A main theme of \cite{lutz2020bounding} is how close this upper bound is to being a lower bound for the points on a line, the answer to which is expressed in the following theorem. Let
\begin{equation*}
    \dim(x \vert y) := \liminf_{r \to \infty} \frac{K_r(x \vert y)}{r}
\end{equation*}
be the \define{conditional dimension of $\bm{x}$ given $\bm{y}$}.

\begin{thm}[N.~Lutz--Stull \cite{lutz2020bounding}] \label{thm:bounding}
    For all $a,b,x \in \R$ and $A \subseteq \{0,1\}^*$,
    \begin{equation*}
        \dim^A(x \vert a,b) + \min \big\{ \dim^A(a,b), \dim^{a,b}(x) \big\} \leq \dim^A(x,ax+b).
    \end{equation*}
\end{thm}

This is the key ingredient in our proof of Proposition \ref{prop:generalized-lse-R2}. However, as the effective packing dimension analogue of this statement is false,\footnote{The authors appreciate Donald Stull informing us of this in private communication.} the proof of Theorem \ref{thm:packing-lse-R2} will require a different strategy.

\section{Line segment extension in the plane} \label{s:plane}

\subsection{The Hausdorff dimension bound for \texorpdfstring{$\bm{s}$}{s}-Hausdorff extensions} We begin this section by using Theorem \ref{thm:bounding} of \cite{lutz2020bounding} to prove a bound on the Hausdorff dimension of line segment extensions in the plane. This proof takes Lutz and Stull's result as a black box and illustrates the connection between effective and classical theorems in this setting, which we will need to prove Theorem \ref{thm:packing-lse-R2}.

\begin{prop-nonum}[\ref{prop:generalized-lse-R2}, Restated]
    Let $s \in (0,1]$ and let $E \subseteq \R^2$ be a union of (at least) $s$-Hausdorff-dimensional subsets of lines. Then $\hdim L_s^{\rm H}(E) \leq \hdim E + 2 - 2s$. In particular, if $E$ is a union of line segments and $\upbold{L}(E)$ is its lineal extension, then $\hdim \upbold{L}(E) = \hdim E$.
\end{prop-nonum}

\textit{Proof.} With $E_\ell := E \cap \ell$, write
\begin{equation*}
    E = \bigcup_{\ell \in \Lambda} E_\ell,
\end{equation*}
where $\Lambda \subseteq \mathcal{A}(2,1)$ is the family of lines $\ell$ such that $\hdim E_\ell \geq s$. By the point-to-set principle,
\begin{equation*}
    \hdim E = \min_{A \subseteq \{0,1\}^*} \sup_{z \in E} \dim^A(z) = \min_{A \subseteq \{0,1\}^*} \sup_{\ell \in \Lambda} \+ \sup_{z \in E_\ell} \dim^A(z)
\end{equation*}
and
\begin{equation*}
    \hdim L_s^{\rm H}(E) = \min_{A \subseteq \{0,1\}^*} \sup_{z \in L_s^{\rm H}(E)} \dim^A(z) = \min_{A \subseteq \{0,1\}^*} \sup_{\ell \in \Lambda} \+ \sup_{z \in \ell} \dim^A(z).
\end{equation*}
Comparing the right-hand sides of these equations, we see it suffices to show that
\begin{equation} \label{eq:effective-hausdorff-lse-R2}
    \sup_{z \in \ell} \dim^A(z) \leq \sup_{z \in E_\ell} \dim^A(z) + 2 - 2s
\end{equation}
for every oracle $A \subseteq \{0,1\}^*$ and every line $\ell \in \Lambda$. Taking such an $A$ and $\ell$, we interchange the $x$- and $y$-coordinates if necessary so that $\ell$ is not vertical and we let $(a,b)$ be the slope-intercept pair of $\ell$. As observed in the previous section, for each $x \in \R$ and each precision $r \in \N$,
\begin{equation*}
    K_r^A(x,ax+b) \leq K_r^A(x,a,b) + O_{x,a,b}(\log r) \leq r + K_r^A(a,b) + O_{x,a,b}(\log r).
\end{equation*}
Hence $\dim^A(x,ax+b) \leq \min \, \{ 1 + \dim^A(a,b), 2 \}$ and, consequently,
\begin{equation} \label{eq:sup-dim(z)-upper-bd}
    \sup_{z \in \ell} \dim^A(z) \leq \min \, \{ 1 + \dim^A(a,b), 2 \}.
\end{equation}
Now, let $S$ be the projection of $E_\ell$ onto the $x$-axis. Then $\hdim S \geq s$, so by the point-to-set principle, for every $\eps>0$, there exists $x_\eps\in S$ such that $\dim^{A,a,b}(x_\eps) \geq s - \eps$. Applying Theorem \ref{thm:bounding}, we have
\begin{align*}
    \dim^A(x_\eps, ax_\eps + b) &\geq \dim^A(x_\eps\vert a, b) + \min \, \{ \dim^A(a,b), \dim^{a,b}(x_\eps) \} \\
    &\geq \dim^{A,a,b}(x_\eps) + \min \, \{ \dim^A(a,b), \dim^{A,a,b}(x_\eps) \} \\
    &\geq s - \eps + \min \, \{ \dim^A(a,b), s-\eps \}.
\end{align*}
Letting $\eps$ go to zero gives
\begin{equation} \label{eq:sup-dim(z)-lower-bd}
    \sup_{z \in E_\ell} \dim^A(z) \geq \min \, \{s + \dim^A(a,b), 2s\}.
\end{equation}
The difference between the upper bound in \eqref{eq:sup-dim(z)-upper-bd} and the lower bound in \eqref{eq:sup-dim(z)-lower-bd} is greatest when $\dim^A(a,b) \geq 1$, so subtracting the latter inequality from the former implies the desired inequality \eqref{eq:effective-hausdorff-lse-R2}. \textqed

\subsection{The packing dimension bound for \texorpdfstring{$\bm{1}$}{1}-Hausdorff extensions} Now we turn our attention to the packing dimension version of the line segment extension problem. The key issue is that we do not have an analogue of Lutz and Stull's bound for effective packing dimension. In fact, the statement obtained by replacing effective Hausdorff dimension with effective packing dimension in Theorem \ref{thm:bounding} is simply not true. 

Essentially, the inequality fails because a high packing dimension for the pair $(a,b)$ can be the result of $K_r^A(a,b)$ growing relatively slowly in $r$ up to some level of precision, and then significantly more quickly up to a higher level of precision. At a key technical step in the proof, the complexity function of $(a,b)$ needs to have certain properties, which can be guaranteed by reducing its complexity up to precision $r$ via an oracle $D$. This ``wastes'' complexity growth of $(x,a,b)$ that we would like to transfer to $(x,ax+b)$, but since effective Hausdorff dimension only reflects a \textit{lower} bound on the asymptotic complexity growth, $D$ does not reduce the complexity of $(x,ax+b)$ unacceptably in comparison to the effective Hausdorff dimension of $(a,b)$. Effective packing dimension, however, reflects an \textit{upper} bound on asymptotic complexity growth, which dashes any hope for the analogous packing inequality. We will proceed without proving an explicit lower bound on the packing dimension of arbitrary points on a line, but will still show that for $x$ of essentially maximal complexity at certain precisions, $(x,ax+b)$ also has essentially maximal complexity. This will imply (a somewhat stronger version of) the line segment extension conjecture for packing dimension in the plane.

\begin{thm-nonum}[\ref{thm:packing-lse-R2}, Restated]
    If $E \subseteq \R^2$ is a union of $1$-Hausdorff-dimensional subsets of lines, then $\pdim L_1^{\rm H}(E) = \pdim E$. In particular, if $E$ is a union of line segments and $\upbold{L}(E)$ is its lineal extension, then $\pdim \upbold{L}(E) = \pdim E$.
\end{thm-nonum}

The proof will proceed in three main steps. First, we will need to understand the complexity function $K_s^A(a,b)$ on the interval $[1,r]$. Based on its behavior, we will judiciously choose an oracle $D$ that reduces the complexity at precisions close to $r$. Next, we will show that with the addition of this oracle, we can apply a technical lemma from \cite{lutz2020bounding} to establish a lower bound for $K_r^A(x,ax+b)$ on  $[1,r]$. This lower bound improves if $x$ has high complexity at all precisions relative to $A,a,b$. Finally, we prove that for such $x$, the lower bound in the previous step essentially matches an upper bound at all sufficiently large precisions, a property which is guaranteed by the choice of $D$.

Before we begin, we will need several lemmas in \cite{lutz2020bounding}, starting with their Lemma 3.1 stated in a relativized form.

\begin{lem} \label{lem:enumeration}
    Suppose $a,b, x \in \R$, $B \subseteq \{0,1\}^*$, $r \in \N$, $\delta \in \R^+$, and $\eps,\eta \in \Q^+$ satisfy $r > \log(2|a| + |x| + 5) + 1$, and the following:
    \begin{enumerate}[label={\normalfont \textbf{(\arabic*)}}, itemsep=1.5pt, topsep=-3pt]
        \item $K_r^B(a,b) \leq (\eta + \eps) r$.
        \item For every $(u,v) \in B((a,b),1)$ such that $ux+v = ax+b$,
        \begin{equation*}
            K_r^B(u,v) \geq (\eta - \eps)r + \delta (r-s)
        \end{equation*}
        whenever $s := -\log |(a,b) - (u,v)| \in (0,r]$.
    \end{enumerate}
    Then 
    \begin{equation*}
        K_r^B(x, ax+b) \geq K_r^B(x,a,b) - \frac{4\eps}{\delta} r - K^B(\eps) - K^B(\eta) - O(\log r).
    \end{equation*}
\end{lem}

We note that the implicit constant may depend on $x$, $a$, and $b$, but these will be fixed in each application. This kind of lemma is often referred to as an ``enumeration'' lemma, as its proof depends on enumerating many short strings to find one that gives an output with the desired properties; enumeration lemmas are key technical elements of many proofs using effective methods because they give us conditions under which a desired lower bound holds. In the proof of our main theorem, showing that the two conditions are satisfied is a significant element in proving the desired lower bound.

We also make use of Lemmas 3.2 and 3.3 from \cite{lutz2020bounding}, stated in the form we will need.\footnote{Lemmas 3.2 and 3.3 are used in a relativized form in \cite{lutz2020bounding}, so we state them in this way. The third and fourth properties in Lemma \ref{lem:finite-oracle}, which are implicit in \cite{lutz2020bounding} Lemma 3.3, are easy consequences of the construction of $D$ and also are enumerated in \cite{stull2022pinned}.} 

\begin{lem} \label{lem:geometric-fact}
    Let $x,a,b \in \R$. For all $(u,v) \in B((a,b),1)$ such that $ux+v = ax+b$ and for all $r \geq s := -\log |(a,b) - (u,v)|$,
    \begin{equation*}
        K_r^A(u,v) \geq K_s^A(a,b) + K_{r-s,r}^A(x \vert a,b) - O(\log r).
    \end{equation*}
\end{lem}

\begin{lem} \label{lem:finite-oracle}
    Let $A \subseteq \{0,1\}^*$, $r \in \N$, $z \in \R^n$, and $\eta \in \Q^+$. There is an oracle $D = D(A, n, r, z, \eta)$ satisfying the following:
    \begin{enumerate}[label={\normalfont \textbf{(\arabic*)}}, itemsep=1.5pt, topsep=-3pt]
        \item For every natural number $t \leq r$, 
        \begin{equation*}
            K_t^{A,D}(z) = \min \, \{ \eta r, K_t^A(z) \} + O(\log r).
        \end{equation*}
        \item For every $m,t \in \N$ and $y \in \R^m$, 
        \begin{equation*}
            K_{t,r}^{A,D}(y \vert z) = K_{t,r}^A(y \vert z) + O(\log r) \quad \text{and} \quad K_t^{A,D,z}(y) = K_t^{A,z}(y) + O(\log r).
        \end{equation*}
        \item If $B \subseteq \{0,1\}^*$ satisfies $K_r^{A,B}(z) \geq K_r^A(z) - O(\log r)$
        \begin{equation*}
            K_r^{A,B,D}(z) \geq K_r^{A,D}(z) - O(\log r).
        \end{equation*}
        \item For every $m,t \in \N$, $u \in \R^n$, and $w \in \R^m$,
        \begin{equation*}
            K_{r,t}^A(u \vert w)\leq K_r^{A,D}(u \vert w) + K_{r,t}^A(z) - \eta r + O(\log r).
        \end{equation*}
    \end{enumerate}
\end{lem}

Lemma \ref{lem:geometric-fact} is the key geometric observation of Lutz and Stull, and it formalizes the statement that lines passing through the same point are either almost parallel (in which case they contain much of the same information), or they are transverse enough their approximations together determine the $x$-coordinate of the intersection to a high precision. 

Lemma \ref{lem:finite-oracle} is common in effective arguments. Although it is lengthy to state, the idea is rather simple: if you want to lower the complexity of a point $z$ at some precision $r$, look back to find a precision $s<r$ at which the complexity of $z$ is what you want it to be at precision $r$. Then, let $D$ encode all of the new information in $z$ from $s$ to $r$. Property 1 says this oracle accomplishes the goal of lowering the complexity. By contrast, the remaining properties tell us that $D$ is not \textit{too} helpful, that is, $D$ does not undesirably lower the complexity of other objects. Specifically $D$ is not any more helpful in any calculation than knowing $z$ up to precision $r$ (Property 2), it does not magically become more useful when combined with other unhelpful oracles for $z$ (Property 3), and it does not reduce the complexity of any object more than it reduces the complexity of $z$ at precision $r$ (Property 4).

\define{Proof of Theorem \ref{thm:packing-lse-R2}.} By the same application of the point-to-set principle used in Proposition \ref{prop:generalized-lse-R2}, it suffices to show that for any planar line $\ell$ with slope-intercept pair $(a,b)$ and for any oracle $A \subseteq \{0,1\}^*$, if $\dim^{A,a,b}(x_\eps) \geq 1-\frac{\eps}{4}$ for a collection of $x_\eps$, then
\begin{equation*}
    \lim_{\eps \to 0} \Dim^A(x_\eps, ax_\eps + b) = \sup_{z \in \ell} \Dim^A(z).
\end{equation*}
With this aim in mind, let $a,b \in \R$ and $A \subseteq \{0,1\}^*$ be given. Let $\eps \in (0,1) \cap \Q$ and assume $x_\eps \in \R$ is such that $\dim^{A,a,b}(x_\eps) \geq 1 - \frac{\eps}{4}$. In the following, we will always assume $r$ is large enough that for $s > \log r$, $(1-\frac{\eps}{2}) s \leq K_s^{A,a,b}(x_\eps)$.

\noindent\textbf{Choosing an oracle for the complexity function on $\bm{[1,r]}$:} For the first step of the argument, our aim is to find, given $A\subseteq\{0,1\}^*$, $a, b\in\R$, and $r$ sufficiently large, an oracle relative to which we can apply Lemma \ref{lem:enumeration} and which does not lower $K_r^A(a, b)$ too much. Let $c_r$ be the largest minimizer of $K_t^A(a, b)-t$ on $[1, r]$. This implies that
\begin{equation*}
    K_s^A(a,b) \geq K_{c_r}^A(a,b) - (c_r-s)  \qquad \forall \+ s \in [1,c_r].
\end{equation*}
This property, which we call the teal property after \cite{stull2022pinned}, is essentially what will allow us to show the second condition of Lemma \ref{lem:enumeration} is satisfied. However, we would like the teal property to hold on $[1,r]$ and not just $[1, c_r]$, which entails reducing $K_r^A(a,b)$ with $D$ so that $K_{r,c_r}^{A,D}(a,b) \approx r - c_r$, i.e., the average growth rate on $[c_r, r]$ is about 1. If possible,\footnote{If no element of the set satisfies the inequalities, just set $\eta=0$. In this case, we can do no better than the trivial lower bound of $0$ in \eqref{eq:D-does-not-decrease} for \textit{this} choice of $\eps$, but in practice, we will pick up any actual growth as we pass through with smaller and smaller $\eps$.} pick $\eta$ to be an element of the finite set $\big\{ \frac{i}{2^m} \!: i \in \N,  m = 2 - \lceil \log \eps \rceil \text{ and } 0 \leq i \leq 2^m \big\}$ such that 
\begin{equation*}
    \dfrac{r - c_r + K_{c_r}^A(a,b)}{r} - 2 \sqrt{\eps} < \eta < \dfrac{r - c_r + K^A_{c_r}(a,b)}{r} - \sqrt{\eps}.
\end{equation*}

\begin{figure}[t]
\centering

\tikzset{every picture/.style={line width=0.75pt}} 

\begin{tikzpicture}[x=0.75pt,y=0.75pt,yscale=-1,xscale=1]

\begin{scope}
    \draw [thick, color={rgb,255: red,0; green,180; blue,180}, -{Stealth}, dash pattern={on 5.5pt off 4.5pt}] (229,290) -- (129,390);
    
    \draw[thick] (125.7,374.3) -- (162.19,290) -- (229,290) -- (288,181.7);

    \fill[black] (229,374.3) circle (2pt);
    \fill[black] (229,290) circle (2pt);
    \fill[black] (288,374.3) circle (2pt);
    \fill[black] (288,181.7) circle (2pt);
    
    \draw[thick,-{Stealth}] (105.7,374.3) -- (304.7,374.3); 
    \draw[thick,-{Stealth}] (125.7,399) -- (125.7,160); 
    
    \draw (172.09,261.61) node [anchor=north west][inner sep=0.75pt] [font=\small] {$K_{s}^{A}(a,b)$};
    \draw (283,378.28) node [anchor=north west][inner sep=0.75pt] [font=\small] {$r$};
    \draw (225.41,378.42) node [anchor=north west][inner sep=0.75pt] [font=\small] {$c_{r}$};
    
\end{scope}

\begin{scope}[shift={(250,0)}]
    \draw [thick, color={rgb,255: red,0; green,180; blue,180}, -{Stealth}, dash pattern={on 5.5pt off 4.5pt}] (288,242.2) -- (129,390);
    
    \draw[thick] (125.7,374.3) -- (162.19,290) -- (229,290) -- (254.8,242.2);
    \draw[thick, color={rgb,150: red,100; green,100; blue,100}] (254.8,242.2) -- (288,181.7); 
    \draw[thick] (254.8,242.2) -- (288,242.2);
    \draw[thick,-{Stealth}] (288,197) -- (288,232);
    
    \fill[black] (229,374.3) circle (2pt);
    \fill[black] (229,290) circle (2pt);
    \fill[black] (288,374.3) circle (2pt);
    \fill[black,color={rgb,150: red,100; green,100; blue,100}] (288,181.7) circle (2pt);
    \fill[black] (288,242.2) circle (2pt);
    
    \draw[thick,-{Stealth}] (105.7,374.3) -- (304.7,374.3); 
    \draw[thick,-{Stealth}] (125.7,399) -- (125.7,160); 

    \draw (168.09,261.61) node [anchor=north west][inner sep=0.75pt] [font=\small] {$K_{s}^{A,D}(a,b)$};
    \draw (283,378.28) node [anchor=north west][inner sep=0.75pt] [font=\small] {$r$};
    \draw (225.41,378.42) node [anchor=north west][inner sep=0.75pt] [font=\small] {$c_{r}$};
    \draw (292,212) node [anchor=west][inner sep=0.75pt] [font=\small] {$D$};
    
\end{scope}

\draw (201.47,402.5) node [anchor=north west][inner sep=0.75pt] [font=\small] {\textbf{(i)}};
\draw (451.47,402.5) node [anchor=north west][inner sep=0.75pt] [font=\small] {\textbf{(ii)}};

\end{tikzpicture}

\caption{\textbf{(i)} The definition of $c_r$ guarantees that the teal property holds on $[1,c_r]$. \textbf{(ii)} We choose $D$ such that the teal property holds on $[1,r]$ relative to $(A,D)$ \textit{and} the average growth rate of $K_s^{A,D}(a,b)$ on $[c_r,r]$ is close to $1$.}
\label{fig:OracleD}
\end{figure}

Using Lemma \ref{lem:finite-oracle}, let $D = D(A,2,r,(a,b), \eta)$. By the definition of $D$, 
\begin{equation*}
    K_r^{A,D}(a,b) = \eta r - O(\log r).
\end{equation*}
Hence, 
\begin{equation} \label{eq:D-does-not-decrease}
    K_r^{A,D}(a,b) \geq K_{c_r}^A(a,b) + (r - c_r) - 2 \sqrt{\eps} r - O(\log r).
\end{equation}
In the next part of the proof, we will use both this lower bound, and the fact that this oracle gives the following version of the teal property:
\begin{equation} \label{eq:better-teal-property}
     K_s^{A,D}(a,b) \geq K_r^{A,D}(a,b) - (1-\sqrt{\eps}) (r-s) - O(\log r).
\end{equation}

To prove \eqref{eq:better-teal-property}, first observe that by definition, 
\begin{equation*}    
    K_s^{A,D}(a,b) = \min \, \big\{ \eta r, K_s^A(a,b) \big\} - O(\log r).
\end{equation*}
In the first case, 
\begin{align*}
    K_s^{A,D}(a,b) &= \eta r - O(\log r) \\
     &= K_r^{A,D}(a,b) - O(\log r) \\
     &\geq K_r^{A,D}(a,b) -  (1-\sqrt{\eps})(r-s)- O(\log r).
\end{align*}
In the second case, by the definition of $c_r$
\begin{equation*}
    K_s^A(a,b) - s \geq K_{c_r}^A(a,b) - c_r.
\end{equation*}
Hence, 
\begin{align*}
    K_s^{A,D}(a,b) &= K_s^A(a,b) - O(\log r)\\
    &\geq K_{c_r}^A(a,b) + s - c_r - O(\log r)\\
    &= K_r^{A,D}(a,b) - \eta r + K_{c_r}^A(a, b) + s - c_r - O(\log r)\\
    &\geq K_r^{A,D}(a,b) - \big( K_{c_r}^A(a,b) + (r - c_r) - \sqrt{\eps} r \big) + K_{c_r}^A(a,b) + s - c_r - O(\log r) \\
    &= K_r^{A,D}(a,b) + s - r + \sqrt{\eps} r - O(\log r) \\
    &\geq K_r^{A,D}(a,b) + s - r + \sqrt{\eps} (r-s) - O(\log r) \\
    &= K_r^{A,D}(a,b) - (1-\sqrt{\eps}) (r-s) - O(\log r)
\end{align*}
and we have established \eqref{eq:better-teal-property}.

\noindent\textbf{Lower bound for $\bm{x_\eps}$:}
Now, we will prove a lower bound on $K_r^A(x_\eps, ax_\eps+b)$ by applying Lemma \ref{lem:enumeration}. If $r$ is sufficiently large, the first condition of Lemma \ref{lem:enumeration} is satisfied for $B=(A, D)$, since 
\begin{equation*}
    K^{A, D}_r(a, b)=\eta r + O(\log r) \leq (\eta +\eps) r.
\end{equation*}
Now we show the second is also satisfied. For $(u,v) \in B\big((a,b), 1\big)$, by Lemma \ref{lem:geometric-fact} and the second property of $D$, we have
\begin{align*}
    K_r^{A,D}(u,v) &\geq K_s^{A,D}(a,b) + K_{r-s,r}^{A,D}(x_\eps \vert a,b) - O(\log r) \\
    &= K_s^{A,D}(a,b) + K_{r-s,r}^A(x_\eps \vert a,b) - O(\log r) \\
    &\geq K_s^{A,D}(a,b) + K_{r-s}^{A,a,b}(x_\eps) - O(\log r).
\end{align*}
Applying \eqref{eq:better-teal-property}, we obtain
\begin{equation*}
    K_{r}^{A,D}(u,v)\geq  K_{r}^{A,D}(a,b) - (1 - \sqrt{\eps})(r-s) + K_{r-s}^{A,a,b}(x_\eps) - O(\log r).
\end{equation*}
Now by our assumption on $x_\eps$, either $r-s\leq \log(r)$ or  $K_{r-s}^{A,a,b}(x_\eps) \geq (1-\frac{\eps}{2}) (r-s)$ holds. In both cases, we have 
\begin{align*}
    K_r^{A,D}(u,v) &\geq K_{r}^{A,D}(a,b) - (1 - \sqrt{\eps})(r-s) + \left( 1 - \frac{\eps}{2} \right) (r-s) - O(\log r) \\
    &= K_r^{A,D}(a,b) - \left( \frac{\eps}{2} - \sqrt{\eps} \right)(r-s) - O(\log r) \\
    &= \eta r - \left( \frac{\eps}{2} - \sqrt{\eps} \right) (r-s)  - O(\log r) \\
    &= \left( \eta - \frac{\eps}{2} + \frac{\eps}{2} \right) r - \left( \frac{\eps}{2} - \sqrt{\eps} \right) (r-s) - O(\log r) \\
    &\geq \left( \eta -\frac{\eps}{2} \right) r - \left( \frac{\eps}{2} - \frac{\eps}{2} - \sqrt{\eps} \right) (r-s) - O(\log r) \\
    &= \left( \eta -\frac{\eps}{2} \right) r + \sqrt{\eps} (r-s) - O(\log r).
\end{align*}
Thus, for sufficiently large $r$, we have
\begin{equation*}
    K_r^{A,D}(u, v)\geq (\eta - \eps) r + \sqrt{\eps}(r-s).
\end{equation*}
This is precisely the second condition of Lemma \ref{lem:enumeration} with $\delta = \sqrt{\eps}$. Both conditions of Lemma \ref{lem:enumeration} are satisfied, hence applying it without any additional oracle, we obtain
\begin{equation*}
    K_r^{A,D}(x_\eps, ax_\eps + b) \geq K_r^{A,D}(x_\eps, a, b) - 4 \sqrt{\eps} r - K(\eps) - K(\eta) - O(\log r).
\end{equation*}
In practice, we will keep the same choice of $\eps$ throughout a partitioning argument even as $r$ goes to infinity, and we chose $\eta$ from a fixed set that depends only on $\eps$. Thus, we can treat the complexity of these terms as constant in $r$. Furthermore, removing an oracle can only increase complexity (up to a log term), so 
\begin{equation*}
    K_r^A(x_\eps, ax_\eps + b) \geq K_r^{A,D}(x_\eps,a,b) - 4 \sqrt{\eps} r - O(\log r).
\end{equation*}
Now, applying symmetry of information and the properties of $D$ to $K_r^{A,D}(x_\eps,a,b)$, we obtain
\begin{align*}
    K_{r}^{A,D}(x_\eps, a, b)&= K_r^{A,D}(x_\eps\vert a, b) + K_r^{A,D}(a,b)  - O(\log r) \\
    &\geq K_r^{A, D, a, b}(x_\eps) + K_r^{A,D}(a,b)  - O(\log r) \\
    &= K_r^{A, a, b}(x_\eps) + K_r^{A,D}(a,b) -  O(\log r) \\
    &\geq \left( 1-\frac{\eps}{2} \right)r + K_r^{A,D}(a,b) -  O(\log r) \\
    &\geq (1-\sqrt{\eps})r + K_r^{A,D}(a,b) -  O(\log r).
\end{align*}
Finally, applying \eqref{eq:D-does-not-decrease}, we establish the desired lower bound when $r$ is sufficiently large:
\begin{equation}\label{eq:new-good-lower-bd}
    K_r^A(x_\eps, ax_\eps + b)  \geq r + K_{c_r}^A(a,b) + (r - c_r)  - 7 \sqrt{\eps} r.
\end{equation}

\noindent\textbf{Upper bound for arbitrary $\bm{x}$:} Now, we want to upper bound $K_r^A(x,ax+b)$ for any given $x\in\R$. We use two facts. The first is that on \emph{any} interval, $K_{r,t}^A(x,ax+b)$ is essentially upper bounded by $2(r-t)$. The second is that $K_r^A(x,ax+b)$ is essentially upper bounded by $K_r^A(x,a,b)$, since precision $r$ approximations of $x, a$, and $b$ are enough to compute $ax + b$ to a similar precision.\footnote{It is \textit{not} true that $K_{r,t}^A(x,ax+b)\leq K_{r,t}^A(x,a,b)$, since $x, a$, and $b$ could all be independently random on $[1,t]$ and then consist only of zeros on $[t,r]$; in this case, the complexity keeps growing for $(x, ax + b)$. Informally this is because the product of $a$ up to precision $t$ and $x$ up to precision $t$ can have length $2t$. This would present a problem for the proof, if we did not have from the definition of $c_r$ that the average growth rate of $K_s^A(a, b)$ is no more than $1$ on $[1,c_r]$ and no less than $1$ on $[c_r,r]$.} Noting that $c_r$ only depends on $A, a, b$, and $r$ (in particular, not on $x$), we want to use the first bound on $[c_r,r]$ and the second on $[1,c_r]$. 

More formally, assume $x$ is such that 
\begin{equation*}
    \Dim^A(x,ax+b) \geq \sup_{z \in \ell} \Dim^A(z) - \eps.
\end{equation*}
Clearly $(x,ax+b) \in \R^2$, so for sufficiently large $r$ and for $t \leq r$, we have
\begin{equation} \label{eq:dimension-partition-upper-bound}
    K_{r,t}^A(x,ax+b) \leq 2(r - t) + o(r-t) \leq 2(r - t) + \eps r.
\end{equation}
At the same time, also for sufficiently large $r$,
\begin{equation} \label{eq:compute-partition-upper-bound}
    K_t^A(x,ax+b) \leq K_t^A(x,a,b) + o(t) \leq K_t^A(x,a,b) + \eps r.
\end{equation}

Now, assume $r$ is large enough that the above bounds hold. We apply \eqref{eq:compute-partition-upper-bound} on the interval $[1, c_r]$ and \eqref{eq:dimension-partition-upper-bound} on $[c_r, r]$. This gives
\begin{align*}
    K_r^A(x,ax+b) &= K_{c_r}^A(x,ax+b) + K_{r,c_r}^A(x,ax+b) + O(\log r) \\
    &\leq K_{c_r}^A(x,a,b) + \eps r + 2(r - c_r) + \eps r + O(\log r) \\
    &\leq c_r + K_{c_r}^A(a,b) + 2(r - c_r) + 3 \eps r.
\end{align*}
Now, assume $\dim^{A,a,b}(x_\eps) > 1 - \frac{\eps}{4}$. By the last step of the proof, this ensures \eqref{eq:new-good-lower-bd} holds. Combining this lower bound with the above upper bound gives that
\begin{equation*}
    K_r^A(x, ax+b) - K_r^A(x_\eps, ax_\eps + b) \leq 10\sqrt{\eps}r
\end{equation*}
for all sufficiently large $r$. Hence, 
\begin{equation*}
    \Dim^A(x,ax+b) - \Dim^A(x_\eps, ax_\eps + b) \leq 10\sqrt{\eps}.
\end{equation*}
Finally, by our choice of $x$, this gives
\begin{equation*}
    \Dim^A(x_\eps, ax_\eps + b) \geq \sup_{z\in\ell} \Dim^A(z) - 11\sqrt{\eps}.
\end{equation*}
We could pick $\eps$ to be arbitrarily small, so this completes the proof. \textqed

Using the lower bound \eqref{eq:new-good-lower-bd}, we also establish (a slightly stronger version of) the generalized Kakeya conjecture for packing dimension in the plane.

\begin{cor-nonum}[\ref{cor:packing-gen-kakeya}, Restated]
    Let $E \subseteq \R^2$ and let $D \subseteq \bbp^1$ be the set of directions of lines intersecting $E$ in a set of Hausdorff dimension $1$. If $D \neq \varnothing$, then
        \begin{equation*}
            \pdim D + 1 \leq \pdim E.
        \end{equation*}
\end{cor-nonum}

\textit{Proof.} Let $\Lambda = \{ \ell \in \mathcal{A}(2,1) \!: \hdim(E \cap \ell) = 1 \}$ and let $F_\ell$ denote the orthogonal projection of $E \cap \ell$ onto the $x$-axis. Without loss of generality, we assume that $\Lambda$ contains no vertical lines. Define $E_1 = \bigcup_\Lambda (E \cap \ell)$. Identifying $D$ with the set of slopes of lines in $\Lambda$, and letting $A$ be a packing oracle for $E_1$ and $D$, we see by the point-to-set principle that 
\begin{equation*}
    \sup_{a \in D} \Dim^A(a) = \pdim D
\end{equation*}
and 
\begin{equation*}
\sup_{\ell \in \Lambda} \sup_{x \in F_\ell} \Dim^A(x,ax+b) = \pdim E_1 \leq \pdim E.
\end{equation*}
Hence, the desired inequality \eqref{eq:packing-gen-kakeya} follows if
\begin{equation*}
    \sup_{a \in D} \Dim^A(a) + 1 \leq \sup_{\ell \in \Lambda} \sup_{x \in F_\ell} \Dim^A(x,ax+b),
\end{equation*}
where the line $\ell$ is given by $y = ax+b$. It therefore suffices to show that for all $a,b \in \R$, $A \subseteq \{0,1\}^*$, and $S \subseteq \R$ of Hausdorff dimension $1$, 
\begin{equation*}
    \Dim^A(a) + 1 \leq \sup_{x\in S} \Dim^A(x,ax+b).
\end{equation*}
By the point-to-set principle, this follows if for any $x_\eps \in \R$ such that $\dim^{A,a,b}(x_\eps) > 1 - \frac{\eps}{4}$,
\begin{equation*}
    \Dim^A(a) + 1 \leq \liminf_{\eps \to 0} \Dim^A(x_\eps, ax_\eps + b).
\end{equation*}
Assuming $r$ is sufficiently large, \eqref{eq:new-good-lower-bd} implies
\begin{align*}
    K_r^A(x_\eps, a x_\eps +b) &\geq c_r + K_{c_r}^A(a,b) + 2(r - c_r) - 7 \sqrt{\eps} r \\
    &\geq c_r + K_{c_r}^A(a) + 2(r - c_r) - 8 \sqrt{\eps} r \\
    &\geq c_r + K_{c_r}^A(a) + (r - c_r) + K_{r,c_r}^A(a) - 9 \sqrt{\eps} r \\
    &\geq c_r + K_r^A(a) + (r - c_r) - 10 \sqrt{\eps} r \\
    &= r + K_r^A(a) - 10 \sqrt{\eps} r.
\end{align*}
This holds at every sufficiently large precision $r$, so dividing by $r$ and taking the limit superior on both sides completes the proof. \textqed

\subsection{Pathological behavior of \texorpdfstring{$\bm{s}$}{s}-packing extensions}

In contrast to $1$-Hausdorff extensions, which do not increase the Hausdorff or packing dimension of $E \subseteq \R^2$, the $1$-packing extensions do not play so well with either notion of dimension. It is clear that the $s$-packing extension can increase the Hausdorff dimension of set: just let $E$ be a Hausdorff dimension $0$, packing dimension-$s$ subset of a line. In fact, the $1$-packing extension can maximally increase the packing dimension of some sets $E$, as illustrated in the following example. 

\begin{eg}[$\bm{L_1^{\rm P}}$ can increase the packing dimension] \label{eg:1-packing-ext} Let $\{r_i\}_{i \in \N}$ be a rapidly increasing sequence. Define 
\begin{align*}
    X &= \big\{ x = 0.x_1x_2x_3 \dots \in [0,1] \!: j \in [r_{2i}, r_{2i+1}) \text{ for some } i \in \N \,\Rightarrow\, x_j = 0 \big\} \quad \text{and} \\
    S &= \big\{a = 0.a_1a_2a_3 \dots \in [0,1] \!: j \in [r_{2i+1}, r_{2i+2}) \text{ for some } i \in \N \,\Rightarrow\, a_j = 0 \big\}.
\end{align*}
Note that $X$ and $S$ both have packing dimension $1$. Further define
\begin{equation*}
    E = \{ (x,ax) \in \R^2 \!: x \in X, \, a \in S \}.
\end{equation*}
Clearly $E$ has packing dimension at least $1$. We get equality using the point-to-set principle: for any $A \subseteq \{0,1\}^*$ relative to which $\{r_i\}_{i \in \N}$ is computable, any $x \in X$, and any $a \in S$, 
\begin{align*}
    K_r^A(x,ax) &\leq K_r^A(x,a) + O(\log r) \\
    &= K_r^A(x+a) + O(\log r) \\
    &\leq r + O(\log r),
\end{align*}
where the second line holds because $x$ and $a$ can only have nonzero digits on disjoint intervals of precision, so their sum is enough to compute both of them. Taking the limit superior shows that the packing dimension of any point in $E$ is no more than $1$, so $\pdim E = 1$.

Since $E$ was defined to be a collection of packing dimension-$1$ subsets of lines with slope $a$, 
\begin{equation*}
    \{ (x,ax) \in \R^2: x \in \R, \, a \in S \} \subseteq L_1^{\rm P}(E).
\end{equation*}
Choose $a \in S$ such that, when the intervals $[r_{2i+1}, r_{2i+2})$ are removed from its binary expansion, the remaining string $\tilde{a}$ is random relative to $A$; since $r_{2i} = o(r_{2i+1})$, this ensures that $\Dim^A(a) = 1$. In addition, let $x$ be random relative to $(A,a)$, so that $\dim^A(x) = \dim^{A,a}(x) = 1$. Because $(x,a)$ can be calculated to precision $r$ from $(x,ax)$ at precision $r$ with little additional information,
\begin{equation*}
     K_r^A(x,ax) = K_r^A(x,a) + o(r) \geq K_r^{A,a}(x) + K_r^A(a) + o(r).
\end{equation*}
Taking the limit superior of both sides as $r \to \infty$ gives
\begin{align*}
    \Dim^A(x,ax) &= \limsup_{r \to \infty} \frac{K_r^A(x,ax)}{r} \\
    &\geq \liminf_{r \to \infty} \frac{K_r^{A,a}(x)}{r} + \limsup_{r \to \infty} \frac{K_r^A(a)}{r} \\
    &= \dim^{A,a}(x) + \Dim^A(a) = 2.
\end{align*}
The oracle $A$ was arbitrary, so $\pdim L_1^{\rm P}(E) = 2$.
\end{eg}

These examples illustrate the lack of nontrivial bounds on the increase of Hausdorff and packing dimension for $1$-packing extensions in the plane. More generally, for $s>0$, we can seek upper bounds on
\begin{equation*}
    \sup_{E \subseteq \R^2} \left( \hdim L_s^{\rm P}(E) - \hdim E \right) \qquad \text{and} \qquad \sup_{E \subseteq \R^2} \left( \pdim L_s^{\rm P}(E) - \pdim E \right).
\end{equation*}
The aforementioned trivial example shows that the first quantity above is at least $1$. Likewise, a slight modification of Example \ref{eg:1-packing-ext} shows that the same is true of the second quantity. However, we would not be surprised by stronger examples in both cases.\footnote{We consider the extreme case of ``two-point'' extensions in the next section.} We note the relationship of this problem to that of dimension estimates on exceptional sets of orthogonal projections. For Hausdorff dimension in the plane, the sharp bound was proved in \cite{ren2023furstenberg} as a consequence of Theorem \ref{thm:ren-wang} \textemdash hence the connection to line segment extension. Analogous to the worse behavior of $s$-packing extensions, bounding the packing dimension of exceptional sets of projections (relative to either $\hdim E$ or $\pdim E$) is a rather different problem; see \cite{jarvenpaa1994upper}, \cite{falconer1996packing}, and \cite{orponen2015exceptional}.

\section{Extension in higher dimensions and related problems} \label{s:higher-dim}

Given a point $z \in \R^n$ on a line, there are two degrees of freedom in choosing collinear points $x,y$ such that $x + t(y-x) = z$ for some $t \in \R$. To prove Proposition \ref{prop:doubling-bound}, we choose $x$ and $y$ such that (along with another condition) $x$'s first coordinate encodes the largest possible quantity of information about $y$. We justify the possibility of such an encoding as a lemma.

\begin{lem} \label{lem:first-digit-encoding}
    For all $y \in \R^n$, $A \subseteq \{0,1\}^*$, and $\eps > 0$, there exists a dense set of points $x \in \R$ such that, for all sufficiently large $r$ (depending on $x$), 
    \begin{equation} \label{eq:first-digit-encoding}
        K_r^A(y \vert x) \leq \max \big\{ K_r^A(y) - (1 - \eps) r, \eps r \big\}.
    \end{equation}
\end{lem}

The idea of the proof is simple: build the point $x \in \R$ such that successive segments of its binary expansion are strings that aid in the computation of successive segments of $y$.

\textit{Proof.} Given a rational $0 < \delta < 1$, we build a point $x_\delta \in \R$ as follows. For each $i \in \N$, let $[r_i,r_{i+1}]$ be an interval of length $\lceil (1+\delta)^i \rceil$, where $r_0 = 1$. Let $\sigma_i$ denote a string testifying to $K_{r_{i+1},r_i}^A(y)$ and let $x_\delta \in [0,1]$ be the real number with binary expansion $0.\sigma_1 \sigma_2 \dots$.

Let $\pi = \pi_1 \pi_2 \pi_3$, where $\pi_1$ is a given string, $U^A(\pi_2)$ is a finite list of positive integer lengths $l_0, \dots, l_k$, and $U^A(\pi_3)$ is a rational number. Define an oracle Turing machine $M^A$ that computes $M^A(\pi)$ as follows. The machine $M^A$ first calculates $U^A(\pi_3) = q \in \Q$ and determines the successive strings $q_0, \dots, q_k$ of lengths $l_0, \dots, l_k$ formed from the binary digits of $q$. It then iteratively computes
\begin{equation*}
    U^A(q_0) = p_0, \ U^A(q_1, p_0) = p_1, \ \dots, \ U^A(q_k, p_{k-1}) = p_k
\end{equation*}
and returns $M^A(\pi) = U^A(\pi_1, p_k)$ as the output. Let $c_M$ be a constant for this machine.

Now let $r \in \N$ be sufficiently large and let $t = r_m$ be the lesser of (1) the largest precision $r_i$ as defined above such that $K_{r_i}^A(y) \leq r$ and (2) the smallest precision $r_i$ such that $r_i \geq r$. (This $r$ can be assumed to be large enough that $K_{r_1}^A(y) \leq r$.) If there is no largest such $r_i$ in (1), we default to the $r_i$ given by (2). As $\{r_i\}_{i\in\N}$ is an increasing sequence tending to infinity, such an $r_i$ clearly exists.

Let $\pi_1$ testify to $K_{r,t}^A(y)$, let $\pi_2$ testify to $K^A\big(K_{r_0}^A(y), K_{r_1,r_0}^A(y), \dots, K_{r_m,r_{m-1}}^A(y) \big)$, and let $\pi_3$ testify to $K_{r'}^A(x_\delta)$, where $r' := K_{r_0}^A(y) + \sum_{i=0}^{m-1} K_{r_{i+1},r_i}^A(y)$. Note that $r' \leq r + O(m \log r)$. One may check that on these inputs, $M^A$ outputs a precision-$r$ estimate of $y$:\footnote{Technically, our approximation may be off by one in its digits. We could resolve this by also giving $M^A$ a short input specifying how to address the ambiguity; this does not increase the input length by more than a logarithmic amount, so for simplicity we ignore the possibility.} from the definition of $x_\delta$, the $q_i$ are precisely the strings that give the additional information in $y$ from precision $r_i$ to precision $r_{i+1}$. From this and by repeatedly using the fact that each $K_{r_{i+1},r_i}^A(y) \in \{ 1, \dots, r' \}$, we find that
\begin{align*}
    K_r^A(y \vert x_\delta) &= K_{r,r'}^A(y \vert x_\delta) + K_{r',r}^A(x_\delta) + O(\log r + \log r') \\
    &\leq K_{r,r'}^A(y \vert x_\delta) + O(m \log r) \\
    &\leq K_{r,t}^A(y) + K^A\big(K_{r_0}^A(y), K_{r_1, r_0}^A(y), \dots, K_{r_m, r_{m-1}}^A(y) \big) + O(m \log r) + c_M \\
    &= K_{r,t}^A(y) + O(m \log r).
\end{align*}
Since $m \leq O(\log r)$ by the definition of $t$, this in turn implies
\begin{equation*}
    K_r^A(y \vert x_\delta) \leq K_{r,t}^A(y) + O((\log r)^2).
\end{equation*}
Now consider the two cases for the choice of $t$. If $t \geq r$, then $K_{r,t}^A(y) \leq O(\log r)$ and for sufficiently large $r$
\begin{equation} \label{eq:high-dimension-lemma-case-1}
    K_r^A(y \vert x_\delta) \leq O((\log r)^2) < 3n \delta r.
\end{equation}
For the case that $t < r$, let $s = r_{m+1}$. Then $K_s^A(y) > r$ and $t \leq s$ by the definition of $t$, so 
\begin{align*}
    K_r^A(y \vert x_\delta) &\leq K_r^A(y) - K_t^A(y) + O((\log r)^2) \\ 
    &= K_r^A(y) - K_s^A(y) + K_{s,t}^A(y)  + O((\log r)^2) \\ 
    &\leq K_r^A(y) - r + K_{s,t}^A(y)  + O((\log r)^2) \\ 
    &\leq K_r^A(y) - r + n(s-t) + O((\log r)^2) 
\end{align*}
by the choices of $t$ and $s$. Since
\begin{equation*}
    s-t = r_{m+1} - r_m \leq 2 \big( (1+\delta)r_m - r_m \big) = 2 \delta r_m < 2 \delta r,
\end{equation*}
it immediately follows that 
\begin{equation*}
    K_r^A(y \vert x_\delta) \leq K_r^A(y) - r + 2n \delta r + O((\log r)^2).
\end{equation*}
So, for all $r$ sufficiently large that $n \delta r$ dominates the $O((\log r)^2)$ term, we have
\begin{equation} \label{eq:high-dimension-lemma-case-2}
      K_r^A(y \vert x_\delta) \leq K_r^A(y) - r + 3n \delta r.
\end{equation}
Picking $\delta < \tfrac{\eps}{3n}$ and combining \eqref{eq:high-dimension-lemma-case-1} with \eqref{eq:high-dimension-lemma-case-2} gives the existence of one $x = x_\delta$ satisfying \eqref{eq:first-digit-encoding}. Appending the digits of $x$ to any dyadic rational (which are dense in $\R$) gives it the same property for sufficiently large $r$, completing the proof. \textqed

For simplicity and considering how it is used below, we stated the lemma for $x \in \R$, but an almost identical proof allows one to encode information about $y \in \R^n$ within $x \in \R^m$.

\begin{prop-nonum}[\ref{prop:doubling-bound}, Restated]
    If $E \subseteq \R^n$ is a union of line segments, then
    \begin{equation*}
        \hdim \upbold{L}(E) \leq \hdim E + \pdim E - 1 \quad \text{and} \quad \pdim \upbold{L}(E) \leq 2 \pdim E - 1.
    \end{equation*}
\end{prop-nonum}

\textit{Proof.} Let $A \subseteq \{0,1\}^*$ be both a Hausdorff oracle and a packing oracle for $E$ and let $\eps > 0$. Given a point $z \in \R^n$, we construct a machine $M^A$ (depending on $z$) that operates as follows. Let $[ \,\cdot\, ]$ denote the fractional part of a real number, and for each $y \in \R^n$ let $c = c_y \in \N$ be such that $[ 2^c y_1 ] = [ 2^c z_1 ]$, i.e., such that the first coordinates of $y$ and $z$ agree from the $c$th binary digit onwards, assuming such a number exists. Finally, given $y \in \R^n$ and given any other $x \in \R^n$ with $x_1 \neq y_1$, we denote
\begin{equation*}
    t = t_{x,y} := \frac{z_1 - x_1}{y_1 - x_1}.
\end{equation*}
Then $z = x + t(y-x)$ with this choice of $t$, provided $x$ lies on the line through $y$ and $z$. Fixing a parameter $s \in \N$ to be specified momentarily, we then let $\pi_1, \dots, \pi_4 \in \{0,1\}^*$ be such that
\begin{enumerate}[label=\textbullet, topsep=-2pt, itemsep=0pt]
    \item $p := U^A(\pi_1)$ is a precision-$(r+s)$ approximation of $x$,
    \item $q := U^A(\pi_2,p)$ is a precision-$(r+s)$ approximation of $y$,
    \item $U^A(\pi_3) = \lceil 2^c z_1 \rceil$ is the first $c$ digits of $z_1$ (equivalently, of $y_1$) for some $c \in \N$ as above, and
    \item $u := U^A\big( \pi_4, p, q, \lceil 2^c z_1 \rceil \big)$ is a precision-$(r+s)$ approximation of $t$.
\end{enumerate}
The machine $M^A$ is taken to satisfy
\begin{equation*}
    M^A(\pi) = p + u(q-p), \quad \text{where} \quad \pi = \pi_1\pi_2\pi_3\pi_4.
\end{equation*}
In order for $M^A(\pi)$ to be a precision-$r$ approximation of $z$, we select $s$ large enough that $p' + u'(q'-p') \in B(z,2^{-r})$ whenever $p',q' \in \Q^n$ and $u' \in \Q$ are precision-$(r+s)$ approximations of $x,y \in \R^n$ and $t \in \R$, respectively. Notice that $s$ depends only on $x$, $y$, and $z$, but in particular not on $r$.

Now let $z \in \ell$ for some line $\ell \subseteq \upbold{L}(E)$ and let $y \in E \cap \ell$ be a point such that, up to a permutation of the axes, $[ 2^c y_1 ] = [ 2^c z_1 ]$ for some $c \in \N$. Next, let $x \in E \cap \ell$ be a point with $x_1 \neq y_1$ such that $x_1$ assists in the computation of $y$ as in Lemma \ref{lem:first-digit-encoding}, i.e., such that
\begin{equation*}
    K_r^A(y \vert x_1) \leq \max \big\{ K_r^A(y) - (1 - \eps) r, \eps r \big\}
\end{equation*}
for all sufficiently large $r$. Up to a loss of $O(\log r)$, the same inequality holds with $x$ in place of $x_1$, the other $n-1$ coordinates of $x$ being ignored in the computation of $y$. Finally, let $t$ be as before, and let $\pi_1$, $\pi_2$, $\pi_3$, and $\pi_4$ witness $K_{r+s}^A(x)$, $K_{r+s}^A(y \vert x)$, $K^A\big( \lceil 2^c z_1 \rceil \big)$, and $K_{r+s}^A\big( t \vert x,y, \lceil 2^c z_1 \rceil \big)$, respectively. Since $t$ is computable from $x$, $y$, and $\lceil 2^c z_1 \rceil$, we have
\begin{equation*}
    K_{r+s}^A\big( t \vert x, y, \lceil 2^c z_1 \rceil \big) = o(r+s) = o(r).
\end{equation*}
Hence, by the design of $M^A$, the choice of $x$, and symmetry of information, the following holds for all large $r \in \N$:
\newpage
\begin{align*}
    K_r^A(z) &\leq |\pi_1 \pi_2 \pi_3 \pi_4| + c_M \\
    &\leq K_{r+s}^A(x) + K_{r+s}^A(y \vert x) + K^A\big( \lceil 2^c z_1 \rceil \big) + K_{r+s}^A\big( t \vert x,y, \lceil 2^c z_1 \rceil \big) + O(\log(r+s)) + c_M \\
    &\leq K_{r+s}^A(x) + K_{r+s}^A(y \vert x) + o(r) \\
    &\leq K_r^A(x) + K_r^A(y \vert x) + 2sn + o(r) \\
    &\leq K_r^A(x) + \max \big\{ K_r^A(y) - (1 - \eps) r, \eps r \big\} + o(r).
\end{align*}
Dividing through by $r$ and taking the limit inferior as $r \to \infty$ gives
\begin{align*}
    \dim^A(z) &\leq \dim^A(x) + \max \big\{ \Dim^A(y) - (1 - \eps), \eps \big\} \\
    &\leq \hdim E + \pdim E - (1 - \eps),
\end{align*}
where we have chosen the first alternative in the maximum because $\pdim E \geq 1$. Taking the supremum over all $z \in F$ gives the first inequality in \eqref{eq:doubling-bound} modulo an $\eps$, which we let decrease to $0$. To obtain the second inequality in \eqref{eq:doubling-bound}, we simply take the limit superior instead of the limit inferior. \textqed

On some level, the argument is morally similar to that of \cite{fraser2017some} Theorem 6, which leverages the dimension inequalities for product sets. Their Kakeya set estimate
\begin{equation*}
    \hdim K + \pdim K \geq n+1
\end{equation*}
improves on Corollary \ref{cor:kakeya-bound} by $+1$, but this is to be expected, as they prove their estimate directly rather than by way of line segment extension.

The complications in the proof of Proposition \ref{prop:doubling-bound} arise primarily from the application of the encoding lemma. Forgoing this encoding, one can one can modify the conclusions of Proposition \ref{prop:doubling-bound} to
\begin{gather*}
    \hdim \upbold{L}(E) \leq \min \big\{ \hdim E + \pdim D, \+ \hdim D + \pdim E \+ \big\} \quad \text{and} \\
    \pdim \upbold{L}(E) \leq \pdim E + \pdim D,
\end{gather*}
where $D \subseteq \bbp^{n-1}$ is the set of directions of the segments in $E$. Specifically, one computes a given $z = x+tv \in \ell \subseteq \upbold{L}(E)$ from $x \in E \cap \ell$, $v \in D$, and $t \in \R$, where $x$ is chosen such that $\Dim(t) = 0$. With this setup, there are no longer enough degrees of freedom to encode additional information about $z$ in one of these parameters as in Lemma \ref{lem:first-digit-encoding}\textemdash hence the disappearance of the $-1$ terms. Bounding the dimension of $D$ in terms of the dimension of $E$ is the generalized Kakeya problem, so in practice these inequalities are no more useful than those in \eqref{eq:doubling-bound}. Similarly, playing with the choices of $x$ and $y$ in the proof of Proposition \ref{prop:doubling-bound} can allow one to derive similar inequalities that only lead to better estimates on the dimension of $\upbold{L}(E)$ given additional structural information about $E$.

In this vein, at the cost of control over the base point $x$ and the scalar $t$, one can turn the inequality for line segment extensions into an inequality for ``two-point extensions" (recall Definition \ref{defn:extensions}).

\begin{prop} \label{prop:2pt-ext}
    If $E \subseteq \R^n$, then
    \begin{equation} \label{eq:2pt-ext}
        \hdim L_0(E) \leq \hdim E + \pdim E + 1 \quad \text{and} \quad \pdim L_0(E) \leq 2 \pdim E + 1.
    \end{equation}
\end{prop}

\textit{Proof.} The proof is the same in spirit as that of Proposition \ref{prop:doubling-bound}, but forgoing the encodings greatly simplifies the matter. Let $z \in L_0(E)$, so that $z = x + t(y-x)$ for some $x,y \in E$ and some $t \in \R$. Taking the limit inferior of both sides of
\begin{equation*}
    K_r^A(z) \leq K_r^A(x) + K_r^A(y) + K_r^A(t) + o(r) \leq K_r^A(x) + K_r^A(y) + r + o(r)
\end{equation*}
gives the first inequality and taking the limit superior gives the second. \textqed

\begin{eg}[Sharpness of Proposition \ref{prop:2pt-ext}] \label{eg:2pt-ext-sharp}
    The inequalities in \eqref{eq:2pt-ext} are sharp in $\R^2$. Let $C_\alpha \subset [0,1]$ be the middle-$\alpha$ Cantor set, $\alpha \in \big( \tfrac{1}{2}, 1 \big)$, so that
    \begin{equation*}
        s := \hdim C_\alpha = \pdim C_\alpha = \frac{\log \tfrac{1}{2}}{\log \tfrac{1}{2}(1 - \alpha)} \in \big( 0, \tfrac{1}{2} \big).
    \end{equation*}
    Then $\mathcal{H}^{2s}(C_\alpha \times C_\alpha) > 0$, so by Marstrand's projection theorem, $\hdim(tC_\alpha - C_\alpha) = 2s$ for $\mathcal{L}^1$-a.e.~$t \in \R$. In fact, by Proposition 1.3 of \cite{peres2000self}, it is also the case that $\pdim(tC_\alpha - C_\alpha) = 2s$ for a.e.~$t \in \R$, so we fix a $t$ satisfying both these equations and let
    \begin{equation*}
        E := (\{0\} \times C_\alpha) \cup (\{1\} \times t C_\alpha).
    \end{equation*}
    Then the set of slopes of lines in $L_0(E)$ is simply $(tC_\alpha - C_\alpha) \cup \{\infty\}$, and in particular the set $D \subset \bbp^1$ of directions in which $L_0(E)$ contains a line has both Hausdorff and packing dimension $2s$. By the generalized Kakeya conjecture in the plane (cf.~Conjecture \ref{conj:gen-kakeya} and Theorem \ref{thm:kakeya-v-gen-kakeya} above),
    \begin{align*}
        \pdim L_0(E) &\geq \hdim L_0(E) \geq \hdim D + 1 = 2s + 1 \\
        &= \hdim E + \pdim E + 1 = 2 \pdim E + 1. 
    \end{align*}
    Hence, both inequalities in \eqref{eq:2pt-ext} hold with equality. Taking an intersection of the exceptional sets in \cite{peres2000self} Proposition 1.3 shows that our argument similarly works when $\alpha = \tfrac{1}{2}$ (although $tC_\alpha - C_\alpha$ may have zero $1$-dimensional packing \textit{measure}), yielding the extreme case $\hdim E = \pdim E = \tfrac{1}{2}$. In the other extreme, any two-point set poses a sharp example for $\hdim E = \pdim E = 0$.
\end{eg}

Interestingly, $\pdim E$ cannot be replaced with $\hdim E$ in the first inequality of Proposition \ref{prop:2pt-ext}. In fact, no nontrivial inequality bounding the Hausdorff dimension of the two-point extension of a set is possible solely in terms of the Hausdorff dimension of the original set.

\begin{eg}[Failure of $\bm{\hdim L_0(E) \leq 2 \hdim E + 1}$] Through a simple argument (in the spirit of \cite{falconer2004fractal} Example 7.8, \cite{hatano1971notes}, or the construction in \cite{altaf2023distance} Theorem 1.4) we observe that the analogous bound for the Hausdorff dimension of two-point extensions severely fails. Let $E = \{ x \in \R^n: \dim(x)=0 \}$. It is immediate by the point-to-set principle that this set has Hausdorff dimension $0$ (although this also follows from a simpler counting argument; see \cite{li2008introduction} Theorem 3.3.1). Its two-point extension is
\begin{align*}
    L_0(E) &= \{ x + t(y-x) \in \R^n \!: \dim(x) = 0, \, \dim(y) = 0, \, t \in \R \} \\
    &\supseteq \{ 2y-x \in \R^n: \dim(x) = 0, \, \dim(y) = 0 \},
\end{align*}
and since scaling a vector by a nonzero rational does not change its pointwise dimension, it follows that
\begin{equation*}
    L_0(E) \supseteq \{ x+y \in \R^n \!: \dim(x) = 0, \, \dim(y) = 0 \}.
\end{equation*}
Now observe that any $z \in \R^n$ can be written as the sum of two Hausdorff dimension-$0$ points. We illustrate for $z \in [0,1]$ with binary representation $0.z_1z_2z_3 \dots$. Let $x = 0.x_1x_2x_3 \dots$, where $x_i = z_i$ when there exists even $j \in \N$ such that $j! \leq i < (j+1)!$ and $x_i = 0$ otherwise. If $y = z - x$, then $x$ and $y$ both consist of alternating blocks of zeros which rapidly increase in length; hence, they both have effective Hausdorff dimension $0$. Repeating the same construction in each coordinate gives the result in $\R^n$. Consequently, $L_0(E) = \R^n$, so
\begin{equation*}
    \hdim L_0(E) = n \quad \text{but} \quad \hdim E = 0.
\end{equation*}
\end{eg}
This also poses a counterexample to any sort of ``reverse continuum Beck's theorem"; see \cite{ren2023discretized} Corollary 1.5.

\phantomsection
\section*{Acknowledgements}
The authors thank Tam\'{a}s Keleti, Donald Stull, Raymond Tana, and Joshua Zahl for helpful comments on an earlier version of this paper. In addition, they extend their gratitude to the referee, whose detailed and insightful comments were of great help during the revision process.

\bibliographystyle{plain}
\bibliography{references}

\end{document}